 \documentclass[12pt]{article}

 \usepackage{mathrsfs,amsfonts,amsmath,amssymb}
 \usepackage[dvips]{color}

 \allowdisplaybreaks

 \newtheorem{theorem}{Theorem}[section]
 \newtheorem{lemma}[theorem]{Lemma}
 \newtheorem{corollary}[theorem]{Corollary}
 \newtheorem{proposition}[theorem]{Proposition}
 \newtheorem{example}[theorem]{Example}

 \def\blemma{\begin{lemma}\sl{}\def\elemma{\end{lemma}}}
 \def\bproposition{\begin{proposition}\sl{}\def\eproposition{\end{proposition}}}
 \def\btheorem{\begin{theorem}\sl{}\def\etheorem{\end{theorem}}}

 \def\beqlb{\begin{eqnarray}}\def\eeqlb{\end{eqnarray}}
 \def\beqnn{\begin{eqnarray*}}\def\eeqnn{\end{eqnarray*}}

 \def\beqlb{\begin{eqnarray}}\def\eeqlb{\end{eqnarray}}
 \def\beqnn{\begin{eqnarray*}}\def\eeqnn{\end{eqnarray*}}

 \def\ar{\!\!\!&}\def\nnm{\nonumber}

 \def\mcr{\mathscr}\def\mbb{\mathbb}\def\mbf{\mathbf}
 \def\mrm{\mathrm}

 \def\proof{\noindent\textit{Proof.~~}}
 
 \def\qed{\quad$\square$\smallskip}
 \definecolor{1}{rgb}{1,0,0}
 \definecolor{1red}{rgb}{0.70,0.00,0.00}

 \def\<{\langle}\def\>{\rangle}

 \def\itDelta{{\it\Delta}}

 \def\d{\mrm{d}}\def\e{\mrm{e}}

 \def\const{\mrm{const}}

\begin{document}

\bigskip\bigskip

\centerline{\Large\bf Fluctuation limits of the super-Brownian}

\medskip

\centerline{\Large\bf motion with a single point catalyst\footnote{
Supported by NSFC grants (No.11126052) and the Fundamental Research
Funds for the Central Universities (No.ZY1116).}}

\bigskip

\centerline{\ Zenghu Li$^a$ ~ and ~ Li
Wang$^{b,}$\footnote{Corresponding author.}}

\medskip

\centerline{\small $^a$Laboratory of Mathematics and Complex
Systems, School of Mathematical Sciences,}

\centerline{\small Beijing Normal University, Beijing 100875, P.R.
China}

\centerline{\small $^b$School of Sciences, Beijing University of
Chemical Technology, Beijing 100029, P.R. China}

\centerline{\small {\tt lizh@bnu.edu.cn} ~ and ~ {\tt
lwang@mail.bnu.edu.cn}}

\bigskip\bigskip

{\narrower{\narrower

\centerline{\bf Abstract}

\bigskip

We prove a fluctuating limit theorem of a sequence of super-Brownian
motions over $\mbb{R}$ with a single point catalyst. The weak
convergence of the processes on the space of Schwarz distributions
is established. The limiting process is an Ornstein-Uhlenbeck type
process solving a Langevin type equation driven by a one-dimensional
Brownian motion.

\bigskip


\bigskip


\par}\par}

\bigskip\bigskip

\section{Introduction}

\setcounter{equation}{0}

In recent years, there has been growing interest in the study of
branching systems in singular media. Although from the viewpoint of
applications some of the models can be artificial, they give useful
insight into the behavior of more realistic systems. The extremely
simple case of the single non-random branching catalyst described by
a Dirac function was introduced in Dawson and Fleischmann (1994).
The model has been studied extensively since then; see, e.g., Dawson
et al. (1995), Fleischmann and Le Gall (1995) and Fleischmann and
Xiong (2006).

In the present paper, we study the fluctuation limits of the single
point catalytic super-Brownian motion (SBM) with small branching.
Our limit theorem shows that the asymptotic fluctuating behavior of
the processes around the Lebesgue measure can be approximated by a
Schwarz distribution-valued Ornstein-Uhlenbeck type process. We also
show that the Ornstein-Uhlenbeck type process solves a Langevin type
equation driven by a one-dimensional Brownian motion. The pathwise
uniqueness for the stochastic equation is established by an explicit
construction of the solution. To prove the weak convergence of the
fluctuating processes on the path space, we first give an extension
of a tightness criterion in Ethier and Kurtz (1986). The results of
this work extend those of Li (2009) and Li and Zhang (2006) on
Dawson-Watanabe superprocesses with immigration; see also Bojdecki
and Gorostiza (1986), Gorostiza and Li (1998) and Dawson et al.
(1989) for some earlier results.

\section{Single point catalytic SBM}

\setcounter{equation}{0}

Let $C(\mbb{R})$ be the Banach space of bounded continuous functions
on $\mbb{R}$ endowed with the supremum norm $\|\cdot\|$. Write
$C_0(\mbb{R})$ for the space of functions in $C(\mbb{R})$ vanishing
at infinity. Let $C^2(\mbb{R})$ denote the space of smooth functions
on $\mbb{R}$ with continuous derivatives up to the second order
belonging to $C(\mbb{R})$. We fix a constant $p>1$ and let
$h_p(x)=(1+x^2)^{-p/2}$ for $x\in \mbb{R}$. Let $C_p(\mbb{R})$
denote the set of continuous functions $f\in C_0(\mbb{R})$
satisfying $|f|\le \const\cdot h_p$ and let $C_p(\mbb{R})^+$ be the
subset of its nonnegative elements. Let $M_p(\mbb{R})$ be the space
of $\sigma$-finite measures $\mu$ on $\mbb{R}$ satisfying
$\int_{\mbb{R}}h_pd\mu< \infty$. Write $\<\mu,f\> = \int_{\mbb{R}}
fd\mu$ for $\mu\in M_p(\mbb{R})$ and $f\in C_p(\mbb{R})$. The
topology on $M_p(\mbb{R})$ is defined by the convention:
 $$
\mu_n\to \mu ~\mbox{if and only if}~\<\mu_n,f\>\to \<\mu,f\>
~\mbox{for all}~ f\in C_p(\mbb{R}).
 $$
We denote the Lebesgue measure on $\mbb{R}$ by $\lambda$, which
clearly belongs to $M_p(\mbb{R})$. Let $(P_t)_{t\ge 0}$ be the
transition semigroup of the one-dimensional standard Brownian motion
$\xi$ generated by $A := \itDelta/2$ and let $\sigma>0$ be a
constant. Let $p(t,x,y) = p(t,y-x)$ denote the transition density of
the Brownian motion. A time-homogeneous Markov process $X = (\Omega,
\mcr{F}, \mcr{F}_t, X_t, \mbf{P}_\mu)$ with state space
$M_p(\mbb{R})$ is called a \textit{SBM with single point catalyst}
at $c\in \mbb{R}$ if it has transition semigroup $(Q_t)_{t\ge 0}$
given by
 \beqlb\label{2.1}
\int_{M_p(\mbb{R})} \e^{-\<\nu,f\>}Q_t(\mu,\d \nu)
 =
\exp\left\{-\<\mu,V_tf\>\right\},
 \eeqlb
where $f \in C_p(\mbb{R})^+$ and $v(t,x) := V_tf(x)$ is the unique
positive solution of the integral evolution equation
 \beqlb\label{2.2}
v(t,x)=P_tf(x)-\frac{\sigma^2}{2}\int_0^tp(t-s,c-x)v(s,c)^2 ds,
\quad t\ge 0, x\in \mbb{R}.
 \eeqlb
The following theorems are generalizations of the results in Dawson
and Fleischmann (1994). In particular, the existence of the single
point catalytic SBM is a consequence of the first theorem.

\btheorem\label{t2.1} The time-homogeneous Markov process $X =
(\Omega, \mcr{F}, \mcr{F}_t, X_t, \mbf{P}_\mu)$ determined by
equation (\ref{2.2}) via the Laplace transition functional
(\ref{2.1}) can be constructed on the space $C([0,\infty),
M_p(\mbb{R}))$ of continuous $M_p(\mbb{R})$-valued trajectories
satisfying $X_t(\{c\})=0$ for all $t>0$. The following expectation
and covariance formulas hold:
 \beqlb
\label{2.3} \mbf{E}_\mu\<X_t,f\>
 \ar=\ar
\<\mu P_t,f\>,\\
\label{2.4} \mbf{Var}_\mu\<X_t,f\>
 \ar=\ar
\sigma^2\int_{\mbb{R}} \mu(dx)\int_0^tp(t-s,c-x)P_sf(c)^2 ds,
 \eeqlb
where $\mu\in M_p(\mbb{R})$ and $f\in C_p(\mbb{R})$. \etheorem

\btheorem\label{t2.2} There is a version of $X$ such that there
exists a jointly continuous random field $x=\{x_t(z): t>0, z\neq
c\}$ satisfying
 $$
X_t(dz)=x_t(z) dz ~~\mbox{for all}~t>0, ~\mbf{P}_\mu\mbox{-a.s.}
 $$
The random field $x$ has the Laplace transforms
 $$
\mbf{E}_\mu\exp\bigg\{-\sum_{i=1}^k x_t(z_i)\theta_i\bigg\}
 =
\exp\left\{-\langle \mu,u(t)\rangle\right\}, \quad t>0, \theta_i\ge
0, z_i\neq c, 1\le i\le k,
 $$
where $u(t,x)\ge 0$ solves
 $$
u(t,x)=\sum_{i=1}^k\theta_ip(t,z_i-x)-\frac{\sigma^2}{2}\int_0^t
p(t-s,c-x)u(s,c)^2 ds, \quad t>0.
 $$
Moreover, the function $f$ in formulas (\ref{2.3}) and (\ref{2.4})
can be replaced by Dirac function $\delta_z$ for any $z\neq c$.
\etheorem

By the sample path continuity of the single point catalytic SBM, we
may introduce the occupation time process $Y=\{Y_t: t\ge 0\}$
related to $X$, defined by
 \beqnn
\langle Y_t, f\rangle=\int_0^t\langle X_s, f\rangle ds, \quad f\in
C_p(\mbb{R})^+.
 \eeqnn
Of course, by the integration, $Y$ is smoother than $X$, and
 \beqlb\label{2.5}
y_t(z):=\int_0^tx_s(z) ds, \quad t\ge 0,z\neq c,
 \eeqlb
yields a density field of $Y$, which is $\mbf{P}_\mu$-a.s. jointly
continuous on $\mbb{R}_+\times\{z\neq c\}$. The next result shows
that we have a everywhere jointly continuous occupation density.

\btheorem\label{t2.3} There is a version of $X$ such that the
density field $y$ of $Y$ defined by (\ref{2.5}) extends continuously
to all of $\mbb{R}_+\times \mbb{R}$. Moreover,
 $$
\mbf{E}_\mu\exp\bigg\{-\sum_{i=1}^k y_t(z_i)\theta_i\bigg\}
 =
\exp\left\{-\langle \mu,u(t)\rangle\right\}, \quad t\ge 0,
\theta_i\ge 0, z_i\in \mbb{R}, 1\le i\le k,
 $$
where $u(t,x)\ge 0$ solves
 $$
u(t,x)=\sum_{i=1}^k\theta_i\int_0^tp(t-s,z_i-x)\ ds -
\frac{\sigma^2}{2}\int_0^tp(t-s,c-x)u(s,c)^2 ds, \quad t\ge 0.
 $$
Moreover, for $s\le t$ and $z\in \mbb{R}$ the following expectation
and variance formulas hold:
 \beqlb
\label{2.6} \mbf{E}_\mu y_t(z)
 \ar=\ar
\int_{\mbb{R}} \mu(dx)\int_0^t p(s, z-x) ds,\\
\label{2.7} \mbf{Var}_\mu y_t(z)
 \ar=\ar
\sigma^2\int_{\mbb{R}}\mu(dx) \int_0^t p(s,c-x)\bigg[\int_s^tp(u-s,
z-c) du\bigg]^2 ds.
 \eeqlb
\etheorem

We call $y_t(z)$ the occupation density of the single point
catalytic SBM at $z\in \mbb{R}$ during the time period $[0,t]$. Set
$\mcr{D}_p(A)=\{f\in C_p(\mbb{R})\cap C^2(\mbb{R}): Af\in
C_p(\mbb{R})\}$.

\btheorem\label{t2.4} For all $f\in\mcr{D}_p(A)$,
 $$
M_t(f):=\langle X_t,f\rangle-\langle X_0,f\rangle-\int_0^t \langle
X_s,A f\rangle ds,\quad t\ge 0,
 $$
is a continuous martingale with quadratic variation process
 $$
\langle M(f)\rangle_t:=\sigma^2 f^2(c)y_t(c), \quad t\ge 0.
 $$
\etheorem

\noindent\textit{Proof of Theorems~\ref{t2.1}.} We here give a
simple construction of the catalytic SBM by summing up an infinite
sequence of processes taking values of finite measures. This
construction is also useful in deriving some properties of the
catalytic SBM. For any $\mu\in M_p(\mbb{R})$ we can find a sequence
of finite measures $\{\mu_i\}_{i\ge 1}$ such that
$\mu=\sum_{i=1}^\infty \mu_i$. For each $i\ge 1$ let $X_i =
\{X_i(t): t\ge 0\}$ be a single point catalytic SBM with initial
measure $\mu_i$. We assume the sequence of processes $X_i$, $i\ge 1$
are defined on the same probability space and are independent. Then
we can define a single point catalytic SBM $X = \{X(t): t\ge 0\}$
with initial measure $\mu$ by
 \beqlb\label{2.8}
X(t) = \sum_{i=1}^\infty X_i(t), \qquad t\ge 0.
 \eeqlb
For $n\ge k\ge 1$ it is easy to see that
 $$
X_{k,n}(t) = \sum_{i=k}^n X_i(t), \qquad t\ge 0
 $$
is a continuous finite measure-valued catalytic SBM with initial
state $\mu_{k,n} := \sum_{i=k}^n\mu_i$. By applying Theorem 1.2.7 in
Dawson and Fleischmann (1994) to the process $X_{k,n}(t)$ we have
 \beqnn
 \ar \ar\mbf{E}\Big[\sup_{0\le s\le t}\<X_{k,n}(s),h_p\>^2\Big]
\le 2\<\mu_{k,n},h_p\>^2 + 16\sigma^2 h_p^2(c)\bigg(\int_{\mbb{R}}
\mu_{k,n}(dx)\int_0^t p(s, c-x) ds \bigg)\cr
 \ar\ar\qquad
+ ~ 4 t\int_0^t\bigg(\sigma^2\int_{\mbb{R}}\mu_{k,n}(dx) \int_0^s
p(s-u,c-x)P_uAh_p(c)^2 du + \<\mu_{k,n} P_s,Ah_p\>^2\bigg)ds.
 \eeqnn
The right hand side tends to zero as $k, n \to \infty$. This implies
that $\{X(t): t\ge 0\}$ can be realized in $C([0,\infty),
M_p(\mbb{R}))$. The moment formulas (\ref{2.3}) and (\ref{2.4})
follow from Theorems~1.2.1 and~1.2.4 in Dawson and Fleischmann
(1994) and the construction (\ref{2.8}). \qed

Now we give three lemmas which will be used in the proof of Theorem
2.2, 2.3 and 2.4. These are modifications of Lemmas~2.6.2, 3.2.1
and~3.2.2 in Dawson and Fleischmann (1994). The proofs are similar
to theirs and are omitted here. Fix $f\in C_p(\mbb{R})$ and set
 $$
u_\theta(t,x) := \theta P_tf(x)-v_\theta(t,x), \quad t\ge 0,
x\in\mbb{R}, \theta\ge 0,
 $$
where $v_\theta(t,x)=V_t(\theta f)(x)$. We denote by
$u_\theta^{(k)}$ the $k$th derivative of $u_\theta$ with respect to
$\theta$ taken at $\theta=0$. For $T\ge 0$ and $f\in C_p(\mbb{R})$
put
 $$
\|Pf\|_T=\sup\{|P_rf(c)|: 0\le r\le T\}.
 $$

\blemma\label{l2.5} For each $k\ge 2$ there is a constant $c_k>0$ so
that the power series $\sum_{k\ge 2}c_k\theta^k$ has a positive
radius of convergence and that
 $$
|u_\theta^{(k)}(t,x)|
 \le
2^{1-k}k! c_k\sigma^{2(k-1)}\|Pf\|^k_T t^{(k-1)/2}\exp
\bigg\{-\frac{|c-x|^2}{2t}\bigg\},
 $$
where $x\in\mbb{R}$, $0\le t\le T$ and $k\ge 2$. \elemma

Set $Z_t = X_t-X_0P_t$. Then we have $\mbf{E}_\mu Z_t = 0$.

\blemma\label{l2.6} For each $k\ge 2$ there exists a constant
$C_k\ge 0$ such that
 $$
|\mbf{E}_\mu\<Z_t, f\>^k|
 \le
C_k t^{k/4}\|Pf\|^k_T\sum_{i=1}^{k-1}\left\<\mu,\exp
\bigg\{-\frac{|c-\cdot|^2}{2t}\bigg\}\right\>^i,
 $$
where $0\le t\le T$, $\mu\in M_p(\mbb{R})$ and $f\in C_p(\mbb{R})$.
\elemma

\blemma\label{l2.7} Let $k\ge 1$, $T>0$ and $\mu\in M_p(\mbb{R})$.
Then there exists a constant $C_k(T,\mu)\ge 0$ such that
 $$
\mbf{E}_\mu\<Z_{t+h}-Z_t, f\>^{2k}
 \le
C_k(T,\mu) \left(\|P(P_hf-f)\|_T^{2k}+h^{k/2}\|Pf\|^k_T\right),
 $$
where $0\le t\le t+h\le T$ and $f\in C_p(\mbb{R})$. \elemma

\medskip

\noindent{\textit{Proof of Theorems~\ref{t2.2}, \ref{t2.3}
and~\ref{t2.4}.}} The existence and the characterizations of the
Laplace transforms of the density fields $x = \{x_t(z): t>0, z\neq
c\}$ and $y = \{y_t(z): t>0, z\in\mbb{R}\}$ follow from the
construction (\ref{2.8}) of the catalytic SBM. The moment formulas
can be derived from the Laplace transforms. The continuity
properties of the fields follow by using Klomogorov's criterion and
the above three lemmas. The martingale problem characterization also
follows from the construction (\ref{2.8}). We leave the details to
the reader. \qed

\section{A tightness criterion}

\setcounter{equation}{0}

A tightness criterion based on the martingale problems was given in
Ethier and Kurtz (1986, p.145). However, the martingales considered
there have absolutely continuous increasing processes. In this
section, we give a generalized version of the result. Although the
proof is similar to that of Ethier and Kurtz (1986), we give it here
for reader's convenience.

Let $E$ be a metric space and let $\bar{C}(E)$ be the space of
bounded and uniformly continuous functions on $E$. For each index
$\alpha$, let $X_\alpha$ be a process with sample paths in
$D([0,\infty),E)$ defined on a probability space
$(\Omega_\alpha,\mcr{F}^\alpha, P_\alpha)$ and adapted to a
filtration $\{\mcr{F}_t^\alpha\}$. Let $A_\alpha$ be an increasing
process which is adapted to the filtration $\{\mcr{F}_t^\alpha\}$
and satisfies
 \beqlb\label{3.1}
\lim_{\delta\rightarrow 0}\sup_\alpha\mbf{E}\Big[\sup_{0\le r\le
T}|A_\alpha(r+\delta)-A_\alpha(r)|\Big] =0.
 \eeqlb
Let $\mcr{L}_\alpha$ denote the Banach space of real-valued
$\{\mcr{F}_t^\alpha\}$-progressively measurable processes with norm
$\|Y\| = \sup_{t\ge 0} \mbf{E}[|Y(t)|]< \infty$. Given $T\ge 0$ and
$h_\alpha\in \mcr{L}_\alpha$, define
 $$
\|h_\alpha\|_{p,T}
 =
\bigg[\int_0^T|h_\alpha(t)|^pdA_\alpha(t)\bigg]^{1/p}
 $$
for $0<p<\infty$ and define $\|h_\alpha\|_{\infty,T} =
{\rm{ess}}\sup_{0\le t\le T}|h_\alpha(t)|$. Let
 $$
\mcr{A}_\alpha=\Big\{(Y,Z)\in \mcr{L}_\alpha\times\mcr{L}_\alpha:
Y(t)-\int_0^t Z(s) d A_\alpha(s) ~\mbox{is an}~
\{\mcr{F}_t^\alpha\}\mbox{-martingale}\Big\}
 $$
Let $Q$ denote the set of rational numbers. Then we have

\btheorem\label{t3.1} Suppose that $C_a$ is a subalgebra of
$\bar{C}(E)$. Let $D$ be the collection of $f\in \bar{C}(E)$ such
that for every $\varepsilon>0$ and $T>0$ there exist
$(Y_\alpha,Z_\alpha)\in \mcr{A}_\alpha$ with
 \beqlb\label{3.2}
\sup_\alpha \mbf{E}\bigg[\sup_{t\in[0,T]\cap Q}|Y_\alpha(t)
-f(X_\alpha(t))|\bigg]<\varepsilon
 \eeqlb
and
 \beqlb\label{3.3}
\sup_\alpha\mbf{E}[\| Z_\alpha\|^p_{p,T}]<\infty ~~ \mbox{for some}
~~ p\in (1,\infty].
 \eeqlb
If $C_a$ is contained in $D$, then $\{f\circ X_\alpha\}$ is tight in
$D([0,\infty),\mbb{R})$ for each $f\in C_a$. \etheorem

\proof Let $\varepsilon>0$ and $T>0$. For $f\in C_a$ we have
$(Y_\alpha,Z_\alpha)\in \mcr{A}_\alpha$ such that (\ref{3.2}) and
(\ref{3.3}) hold. Since $f^2\in C_a$, there are
$(Y'_\alpha,Z'_\alpha)\in \mcr{A}_\alpha$ such that
 $$
\sup_\alpha \mbf{E}\bigg[\sup_{t\in[0,T+1]\cap Q}|Y'_\alpha(t) -
f^2(X_\alpha(t))|\bigg]< \varepsilon
 $$
and
 $$
\sup_\alpha\mbf{E}[\| Z'_\alpha\|^{p'}_{p',T}]<\infty ~~\mbox{for
some}~~p'\in (1,\infty].
 $$
Let $0<\delta<1$. For each $t\in[0,T]\cap Q$ and $u\in[0,\delta]\cap
Q$ we have
 \beqnn
\ar\ar\mbf{E}\Big[(f(X_\alpha(t+u))-f(X_\alpha(t)))^2\Big|
\mcr{F}_t^\alpha\Big]\nnm\cr
 \ar\ar\qquad
=\mbf{E}\Big[f(X_\alpha(t+u))^2-f(X_\alpha(t))^2\Big|
\mcr{F}_t^\alpha\Big]\nnm\cr
 \ar\ar\qquad\qquad
-~ 2f(X_\alpha(t))\mbf{E}\Big[f(X_\alpha(t+u))-f(X_\alpha(t))\Big|
\mcr{F}_t^\alpha\Big]\nnm\cr
 \ar\ar\qquad
\le 2\mbf{E}\bigg[\sup_{t\in[0,T+1]\cap Q}|f(X_\alpha(t))^2 -
Y'_\alpha(t)|\bigg|\mcr{F}_t^\alpha\bigg]\nnm\cr
 \ar\ar\qquad\qquad
+~4\|f\|\mbf{E}\bigg[\sup_{t\in[0,T+1]\cap Q}|f(X_\alpha(t)) -
Y_\alpha(t)|\bigg|\mcr{F}_t^\alpha\bigg]\nnm\cr
 \ar\ar\qquad\qquad
+~\mbf{E}\bigg[\sup_{0\le t\le T}\int_t^{t+\delta}|Z'_\alpha(s)|
dA_\alpha(s)\bigg|\mcr{F}_t^\alpha\bigg]\nnm\cr
 \ar\ar\qquad\qquad
+~2\|f\|\mbf{E}\bigg[\sup_{0\le t\le T}\int_t^{t+\delta}
|Z_\alpha(s)|dA_\alpha(s) \bigg|\mcr{F}_t^\alpha\bigg].
 \eeqnn
It follows that
 \beqlb\label{3.4}
\mbf{E}\Big[(f(X_\alpha(t+u))-f(X_\alpha(t)))^2\Big|
\mcr{F}_t^\alpha\Big]
 \le
\mbf{E}\big[\gamma_\alpha(\delta)|\mcr{F}_t^\alpha\big].
 \eeqlb
where
 \beqlb\label{3.5}
\gamma_\alpha(\delta)
 \ar=\ar
2 \sup_{t\in[0,T+1]\cap Q}|f(X_\alpha(t))^2-Y'_\alpha(t)| +
4\|f\|\sup_{t\in[0,T+1]\cap Q}|f(X_\alpha(t))-Y_\alpha(t)|\cr
 \ar \ar
+~\sup_{0\le t\le T}\int_t^{t+\delta}|Z'_\alpha(s)| dA_\alpha(s) +
2\|f\|\sup_{0\le t\le T}\int_t^{t+\delta} |Z_\alpha(s)|
dA_\alpha(s).
 \eeqlb
Note that the inequality (\ref{3.4}) actually holds for all $0\le
t\le T$ and $0\le u\le \delta$ by the right continuity of
$X_\alpha$. Let $1/p+1/q=1$ and $1/p'+1/q'=1$. By (\ref{3.5}) and
H\"older's inequality we have
 \beqnn
\sup_\alpha\mbf{E}[\gamma_\alpha(\delta)]
 \ar\le\ar
2(1+2\|f\|)\varepsilon + B(\delta,T)^{\frac{1}{q'}}
\sup_\alpha\mbf{E}^{\frac{1}{p'}}
\left[\|Z'_\alpha\|^{p'}_{p',T+1}\right]\cr
 \ar \ar
+~2\|f\|B(\delta,T)^{\frac{1}{q}} \sup_\alpha\mbf{E}^{\frac{1}{p}}
\left[\|Z_\alpha\|^p_{p,T+1}\right],
 \eeqnn
where
 \beqnn
B(\delta,T) = \sup_\alpha\mbf{E}\bigg[\sup_{0\le t\le T}
|A_\alpha(t+\delta)-A_\alpha(t)|\bigg].
 \eeqnn
Then we may select $\varepsilon$ depending on $\delta$ in such a way
that
 $$
\lim_{\delta\to 0}\sup_\alpha\mbf{E}[\gamma_\alpha(\delta)]=0.
 $$
Therefore, $\{f\circ X_\alpha\}$ is tight in $D([0,\infty),\mbb{R})$
by Theorem~8.6 in Ethier and Kurtz (1986, pp.137-138). \qed

\section{A fluctuation limit theorem}

\setcounter{equation}{0}

For each integer $k\ge1$, let $\{X_k(t): t\ge 0\}$ be the single
point catalytic super-Brownian motion characterized by (\ref{2.1})
and (\ref{2.2}) with $\sigma^2/2$ replaced by $\sigma^2/2k^2$. Let
$\{x^k_t(z): t>0, z\neq c\}$ and $\{y^k_t(z): t>0, z\in\mbb{R}\}$ be
the corresponding density and occupation density fields. For
simplicity, we assume $X_k(0) =\lambda$, so
 \beqnn
\mbf{E}\<X_k(t), f\> = \<\lambda P_t,f\>=\<\lambda,f\>, \quad t\ge
0, f\in C_p(\mbb{R}).
 \eeqnn
We define a centered signed-measure-valued Markov
process $\{Z_k(t): t\ge 0\}$ by
 \beqnn
Z_k(t) := k(X_k(t)-\lambda), \qquad t\ge 0.
 \eeqnn
Then $\mbf{E}\<Z_k(t), f\>=0$ for each $f\in C_p(\mbb{R})$.

Let $C^\infty(\mbb{R})$ be the set of bounded infinitely
differentiable functions on $\mbb{R}$ with bounded derivatives. Let
$\mcr{S}(\mbb{R})\subset C^\infty(\mbb{R})$ denote the Schwartz
space of rapidly decreasing functions. That is, a function $f\in
\mcr{S}(\mbb{R})$ is infinitely differentiable and for every $k\ge
0$ and every $n\ge 0$ we have
 \beqnn
\lim_{|x|\to \infty}|x|^n\Big|\frac{d^k}{dx^k} f(x)\Big|=0.
 \eeqnn
The topology of $\mcr{S}(\mbb{R})$ is defined by the increasing
sequence of norms $\{p_0, p_1, p_2, \cdots\}$ given by
 \beqnn
p_n(f) = \sum_{0\le k\le n}\sup_{x\in\mbb{R}} (1+|x|^2)^{n/2}
\Big|\frac{d^k}{dx^k} f(x)\Big|.
 \eeqnn
Let $\mcr{S}^\prime(\mbb{R})$ be the dual space of
$\mcr{S}(\mbb{R})$ endowed with the strong topology. Then both
$\mcr{S}(\mbb{R})$ and $\mcr{S}^\prime(\mbb{R})$ are nuclear spaces;
see, e.g., Treves (1967, p.514 and p.530). We can regard $\{Z_k(t):
t\ge 0\}$ as a process taking values from $\mcr{S}^\prime(\mbb{R})$.

\blemma\label{l4.1} For any $t\ge 0$ and $f\in C_p(\mbb{R})$ we have
 \beqlb\label{4.1}
\mbf{E}[\<Z_k(t),f\>^2]
 =
\sigma^2\int_{\mbb{R}} \lambda(dx)\int_0^tp(t-s,c-x)P_sf(c)^2 ds.
 \eeqlb
\elemma

\proof By Theorem~\ref{t2.1}, we have
 \beqnn
\mbf{E}[\<Z_k(t),f\>^2] =\mbf{Var}\<X_k(t),kf\>
 =
\sigma^2\int_{\mbb{R}} \lambda(dx)\int_0^tp(t-s,c-x)P_sf(c)^2 ds
 .\eeqnn
\qed

\blemma\label{l4.2} For any $t\ge 0$ and $f\in \mcr{S}(\mbb{R})$ we
have
 \beqnn
\sup_{k\ge 1}\mbf{E}\Big[\sup_{0\le s\le t}\<Z_k(s),f\>^2\Big]
 \ar\le\ar
 8\sigma^2f(c)^2t\cr \ar\ar
+ ~2t\sigma^2\int_0^tds\int_{\mbb{R}}\lambda(dx)
\int_0^sp(s-u,c-x)P_uAf(c)^2 du.
 \eeqnn
\elemma

\proof By Theorem~\ref{t2.4} we have
 \beqnn
\<X_k(t),f\> = \<\lambda,f\> + M_k(t,f) + \int_0^t \<X_k(s),Af\>
ds,
 \eeqnn
where $\{M_k(t,f): t\ge 0\}$ is a continuous martingale with
increasing process
 \beqlb\label{4.2}
\<M_k(f)\>_t
 =
\frac{\sigma^2}{k^2}f(c)^2y^k_t(c).
 \eeqlb
It is easy to show that for any $f\in \mcr{S}(\mbb{R})$, $Af\in
\mcr{S}(\mbb{R})$ and $\<\lambda, A f\> = 0$. Then we get
 \beqlb\label{4.3}
\<Z_k(t),f\> = kM_k(t,f) + \int_0^t \<Z_k(s),Af\> ds.
 \eeqlb
By Doob's inequality,
 \beqnn
\ar\ar\mbf{E}\Big[\sup_{0\le s\le t}\<Z_k(s),f\>^2\Big] \cr
 \ar\ar\quad
\le 2k^2\mbf{E}\bigg[\sup_{0\le s\le t}|M_k(s,f)|^2\bigg] +
2\mbf{E}\bigg[\bigg(\int_0^t |\<Z_k(s), Af\>| ds\bigg)^2\bigg] \cr
 \ar\ar\quad
\le 8\sigma^2\mbf{E}\bigg[\int_0^t f(c)^2 d y^k_s(c)\bigg] +
2t\int_0^t \mbf{E}\big[\<Z_k(s),Af\>^2\big] ds \cr
 \ar\ar\quad
\le 8\sigma^2f(c)^2\mbf{E}[y^k_t(c)] +
2t\sigma^2\int_0^tds\int_{\mbb{R}}\lambda(dx)
\int_0^sp(s-u,c-x)P_uAf(c)^2 du\cr
 \ar\ar\quad
= 8\sigma^2f(c)^2t+2t\sigma^2\int_0^tds\int_{\mbb{R}}\lambda(dx)
\int_0^sp(s-u,c-x)P_uAf(c)^2 du.
 \eeqnn
That gives the desired estimate. \qed

\blemma\label{l4.3} For any $G\in C^2(\mbb{R})$ and $f\in \mathcal
{D}_p(A)$ we have
 \beqnn
G(\<Z_k(t),f\>)
 \ar=\ar
\int_0^tG^\prime(\<Z_k(s),f\>)\<Z_k(s),Af\> ds\cr
 \ar \ar
+~ \frac{\sigma^2}{2} \int_0^tG^{\prime\prime}(\<Z_k(s),f\>)f(c)^2d
y_s^k(c) + \mbox{mart.}
 \eeqnn
\elemma

\proof By (\ref{4.2}), (\ref{4.3}) and It\^{o}'s formula, it is easy
to see that
 \beqnn
G(\<Z_k(t),f\>)
 \ar=\ar
\int_0^tG^\prime(\<Z_k(s),f\>)\<Z_k(s),Af\> ds\cr
 \ar \ar
+~\frac{\sigma^2}{2}\int_0^tG^{\prime\prime}(\<Z_k(s),f\>)f(c)^2d
y_s^k(c) + \mbox{local mart.}
 \eeqnn
Since the local martingale in the above equality is actually a
square-integrable martingale, we obtain the desired equality.\qed

\blemma\label{l4.4} As $k\to \infty$, $\{y^k_t(c): t\ge 0\}$
converges in distribution on $C([0,\infty),\mbb{R}_+)$ to $\{t:t\ge
0\}$. \elemma

\proof From the moment formula (\ref{2.6}) we have
 \beqlb\label{4.4}
\mbf{E}[y^k_t(c)]=\int\lambda(dx)\int_0^t p(s, z-x) ds=t, \quad t\ge
0.
 \eeqlb
Using (\ref{2.7}) it is not hard to show that
 \beqlb\label{4.5}
\mbf{Var}
[y^k_t(z)]=\frac{\sigma^2}{k^2}\int\lambda(dx)\int_0^tp(s,c-x)
\bigg[\int_s^tp(u, z-c)du\bigg]^2 ds.
 \eeqlb
Then $y^k_t(c)$ converges in probability to $t$ for each fixed
$t>0$. Consequently, $\{y^k_t(c): t\ge 0\}$ converges in finite
dimensional distributions to deterministic process $\{t: t\ge 0\}$
as $k\to \infty$. Further, by similar calculations as in Dawson and
Fleischmann (1994, p.33), we can show that, for any $T>0$
 $$
\mbf{E}[|y_{t+h}^k(c)-y_t^k(c)|^{2n}]\le C_n(T,\lambda) h^{n/2},
\quad 0<t\le t+h\le T,
 $$
where $C_n(T,\lambda)>0$ is a constant depending only on $T$, $n$
and $\lambda$. Therefore, $\{y^k_t(c):t\ge 0\}_{k\ge 1}$ is tight
and $\{y^k_t(c): t\ge 0\}$ converges weakly to $\{t: t\ge 0\}$ as
$k\to \infty$. \qed

\blemma\label{l4.5} For any $T>0$,
 $$
\lim_{\delta\rightarrow 0}\sup_{k\ge 1}\mbf{E}\bigg[\sup_{0\le r\le
T}\left(y^k_{r+\delta}(c)-y_r^k(c)\right)\bigg]=0.
 $$
\elemma

\proof For each $k\ge 1$, since $y_t^k(c)$ is continuous in $t$ and
$\sup_{0\le r\le T}\left(y^k_{r+\delta}(c)-y_r^k(c)\right)$ is
increasing in $\delta$, we have
 $$
\lim_{\delta\rightarrow 0}\sup_{0\le r\le
T}\left(y^k_{r+\delta}(c)-y_r^k(c)\right)=0,
~~\mbf{P}_\lambda\mbox{-a.s.}
 $$
We may assume $0<\delta<1$, then
 \beqnn
\sup_{0\le r\le T}\left(y^k_{r+\delta}(c)-y_r^k(c)\right)\le
y^k_{T+1}(c), ~~\mbf{P}_\lambda\mbox{-a.s.}
 \eeqnn
In view of (\ref{4.4}) we can use dominated convergence theorem to
obtain
 \beqlb\label{4.6}
\lim_{\delta\rightarrow 0}\mbf{E}\bigg[\sup_{0\le r\le T}
\left(y^k_{r+\delta}(c)-y_r^k(c)\right)\bigg]=0.
 \eeqlb
Observe that for each fixed $T>0$,
 $$
\sup_{0\le r\le T}|y^k_{r}(c)-r|\le y^k_{T}(c)+T
 $$
and the family $\{y^k_{T}(c)\}_{k\ge 1}$ is uniformly integrable by
(\ref{4.4}) and (\ref{4.5}). Then $\{\sup_{0\le r\le T}
|y^k_{r}(c)-r|\}_{k\ge 1}$ is uniformly integrable. On the other
hand, by Lemma~\ref{l4.4} we have
 $$
\sup_{0\le r\le T}|y_r^k(c)-r|\to 0 ~\mbox{in probability as}~ k\to
\infty.
 $$
It follows that
 \beqnn
\lim_{k\to \infty}\mbf{E}\bigg[\sup_{0\le r\le T}
|y^k_r(c)-r|\bigg]=0.
 \eeqnn
Then for any given $\varepsilon>0$, there exists $K = K(\varepsilon)
\ge 1$ so that when $k\ge K$,
 \beqnn
\mbf{E}\bigg[\sup_{0\le r\le T+1}\left|y^k_{r}(c)-r\right|\bigg]
 <
\frac{\varepsilon}{4}.
 \eeqnn
Thus for $k\ge K$ and $0<\delta<\varepsilon/2$ we have
 \beqlb\label{4.7}
\ar\ar\mbf{E}\bigg[\sup_{0\le r\le T}
\left(y^k_{r+\delta}(c)-y_r^k(c)\right)\bigg] \cr
 \ar\ar\quad
\le\mbf{E}\bigg[\sup_{0\le r\le T}\left|y^k_{r+\delta}(c) -
(r+\delta)\right|\bigg] + \mbf{E}\bigg[\sup_{0\le r\le T}
\left|y^k_{r}(c)-r\right|\bigg] + \delta\cr
 \ar\ar\quad
\le2\mbf{E}\bigg[\sup_{0\le r\le T+1}\left|y^k_{r}(c) -
r\right|\bigg] + \frac{\varepsilon}{2}
 \le
\varepsilon.
 \eeqlb
By (\ref{4.6}) we can choose $0<\delta_0 = \delta_0(K)<
\varepsilon/2$ so that when $0<\delta<\delta_0$,
 \beqnn
\sup_{1\le k\le K}\mbf{E}\bigg[\sup_{0\le r\le T}
\left(y^k_{r+\delta}(c)-y_r^k(c)\right)\bigg] < \varepsilon.
 \eeqnn
A combination of this and (\ref{4.7}) completes the proof. \qed

\blemma\label{l4.6} The sequence $\{Z_k(t): t\ge 0\}_{k\ge 1}$ is
tight in $C([0,\infty), \mcr{S}^\prime(\mbb{R}))$. \elemma

\proof We shall prove that for every $f\in \mcr{S}(\mbb{R})$ the
sequence $\{\<Z_k(t),f\>: t\ge 0\}_{k\ge 1}$ is tight in
$C([0,\infty),\mbb{R})$, so the result follows by a theorem of
Mitoma (1983); see also Kallianpur and Xiong (1995, p.82). Let
$A_k(t) = t + \frac{\sigma^2}{2}y_t^k(c)$. By Lemma~\ref{l4.3}, for
any $G\in C^\infty(\mbb{R})$ we have
 \beqnn
G(\<Z_k(t),f\>)
 \ar=\ar
\int_0^tG^\prime(\<Z_k(s),f\>)\<Z_k(s),Af\> ds\cr
 \ar \ar
+~\frac{\sigma^2}{2}\int_0^tG^{\prime\prime}(\<Z_k(s),f\>)f(c)^2d
y_s^k(c) + \mbox{mart.} \cr
 \ar=\ar
\int_0^t \bigg[G^\prime(\<Z_k(s),f\>)\<Z_k(s),Af\>b_k(s) \cr
 \ar \ar\quad
+~ G^{\prime\prime}(\<Z_k(s),f\>)f(c)^2h_k(s)\bigg]dA_k(s) +
\mbox{mart.}
 \eeqnn
where $b_k(s)$ and $h_k(s)$ denote the densities of $ds$ and
$\frac{\sigma^2}{2}dy_s^k(c)$ with respect to $dA_k(s)$,
respectively. From (\ref{4.4}) it is elementary to see that
 \beqnn
\ar\ar\sup_{k\ge 1}\mbf{E}\bigg[\int_0^t\Big|G^\prime
(\<Z_k(s),f\>)\<Z_k(s),Af\>b_k(s) \cr
 \ar\ar\qquad\qquad
+ ~ G^{\prime\prime}(\<Z_k(s),f\>) f(c)^2h_k(s)\Big|^2
dA_k(s)\bigg]< \infty.
 \eeqnn
By Theorem~\ref{t3.1} and Lemma~\ref{l4.5}, we infer that
$\{G(\<Z_k(t),f\>): t\ge 0\}_{k\ge 1}$ is tight in
$D([0,\infty),\mbb{R})$. By Lemma~\ref{l4.2} and Chebyshev's
inequality we have
 \beqnn
\sup_{k\ge 1}\mbf{P}\Big[\sup_{0\le s\le t}|\<Z_k(s),f\>| \ge
\alpha\Big] \to 0
 \eeqnn
as $\alpha \to \infty$. Then $\{\<Z_k(t),f\>: t\ge 0\}_{k\ge 1}$
satisfies the compact containment condition and hence it is tight in
$D([0,\infty),\mbb{R})$ by Theorem~9.1 in Ethier and Kurtz (1986,
p.142). Further, $\{\<Z_k(t),f\>: t\ge 0\}\in C([0,\infty),\mbb{R})$
for each $k\geq 1$, $\{\<Z_k(t),f\>: t\ge 0\}_{k\ge 1}$ is tight in
$C([0,\infty),\mbb{R})$ since convergence in the Skorohod topology
is equivalent to locally uniform convergence in
$C([0,\infty),\mbb{R})$.
 \qed

\blemma\label{l4.7} Let $\{Z_0(t): t\ge 0\}$ be any limit point of
$\{Z_k(t): t\ge 0\}$ in the sense of distributions on $C([0,\infty),
\mcr{S}^\prime(\mbb{R}))$. Then for $G\in C^\infty(\mbb{R})$ and
$f\in \mcr{S}(\mbb{R})$ we have
 \beqnn
G(\<Z_0(t),f\>)
 \ar=\ar
\int_0^tG^\prime(\<Z_0(s),f\>)\<Z_0(s),Af\> d s \cr
 \ar \ar
+~\frac{\sigma^2}{2}\int_0^tG^{\prime\prime}(\<Z_0(s),f\>)f(c)^2 d s
+ \mbox{mart.}
 \eeqnn
\elemma

\proof By passing to a subsequence and using the Skorokhod
representation, we may assume $\{Z_k(t): t\ge 0\}$ and $\{Z_0(t):
t\ge 0\}$ are defined on the same probability space and $\{Z_k(t):
t\ge 0\}$ converges a.s.\ to $\{Z_0(t): t\ge 0\}$ in the topology of
$C([0,\infty), \mathscr{S}^\prime(\mbb{R}))$. From Lemma~\ref{l4.3}
we have
 \beqlb\label{4.8}
G(\<Z_k(t),f\>)
 \ar=\ar
\int_0^tG^\prime(\<Z_k(s),f\>)\<Z_k(s),Af\> d s
+\frac{\sigma^2}{2}\int_0^t G^{\prime\prime}(\<Z_k(s),f\>)f(c)^2 d
s\cr
 \ar \ar
+
~\frac{\sigma^2}{2}f(c)^2\bigg(\int_0^tG^{\prime\prime}(\<Z_k(s),f\>)d
y_s^k(c)-\int_0^tG^{\prime\prime}(\<Z_k(s),f\>) d s\bigg)\nnm\\
 \ar \ar
+~ \mbox{mart.}
 \eeqlb
Let $0=s_0<s_1<\cdots<s_{n-1}<s_n=t$ be a partition of $[0,t]$ so
that $\max_{1\le i\le n}(s_i-s_{i-1})\to 0$ as $n\to\infty$. Then we
have
 \beqlb\label{4.9}
\ar \ar\mbf{E}\bigg[\bigg|\int_0^tG^{\prime\prime}(\<Z_k(s),f\>)d
y_s^k(c)-\int_0^tG^{\prime\prime}(\<Z_k(s),f\>)d s\bigg|\bigg]\cr
 \ar \ar\quad
=\mbf{E}\bigg[\bigg|\lim_{n\to\infty}\sum_{i=1}^n
G^{\prime\prime}(\<Z_k(s_i),f\>) \Big((y_{s_i}^k(c) -
y_{s_{i-1}}^k(c)) - (s_i-s_{i-1})\Big) \bigg|\bigg]\cr
 \ar \ar\quad
\le \|G^{\prime\prime}\|\mbf{E}\bigg[\liminf_{n\to\infty}
\sum_{i=1}^n \bigg|(y_{s_i}^k(c)-y_{s_{i-1}}^k(c))
-(s_i-s_{i-1})\bigg|\bigg]\cr
 \ar \ar\quad
\le\|G^{\prime\prime}\|\liminf_{n\to\infty}\sum_{i=1}^n \mbf{E}
\bigg[\bigg|(y_{s_i}^k(c)-y_{s_{i-1}}^k(c))-(s_i-s_{i-1})
\bigg|\bigg]\cr
 \ar \ar\quad
\le\|G^{\prime\prime}\|\liminf_{n\to\infty}\sum_{i=1}^n
\big[\mbf{Var}(y_{s_i}^k(c)-y_{s_{i-1}}^k(c))\big]^{1/2}.
 \eeqlb
Using (\ref{2.7}) it is not hard to show that
 \beqnn
\mbf{Var}\big(y_{s_i}^k(c)-y_{s_{i-1}}^k(c)\big)
 \ar=\ar
\frac{\sigma^2}{k^2}\int_{\mbb{R}}\lambda(dx)\int^{s_i}_{s_{i-1}}
p(s,c-x) \bigg[\int_s^{s_i}p(u-s, c-c)du\bigg]^2ds \cr
 \ar\le\ar
\frac{2\sigma^2}{k^2\pi}(s_i-s_{i-1})^2.
 \eeqnn
Then the right hand side of (\ref{4.9}) tends to zero as
$k\to\infty$. From (\ref{4.1}) it is easy to show that for any $f\in
\mcr{S}(\mbb{R})$, the sequence $\{\<Z_k(s),f\>\}_{k\ge 1}$ is
uniformly integrable on $ \Omega\times[0, t]$ relative to the
product measure $\mbf{P}(d\omega)ds$. Letting $k\rightarrow\infty$
in (\ref{4.8}) we obtain the desired result.\qed

\bproposition\label{p4.8} For every $\mu\in \mcr{S}^\prime(\mbb{R})$
there is a process $\{Z(t): t\ge 0\}$ with sample paths in
$C([0,\infty), \mcr{S}^\prime(\mbb{R}))$ so that for $G\in
C^\infty(\mbb{R})$ and $f\in \mcr{S}(\mbb{R})$ we have
 \beqlb\label{4.10}
G(\<Z(t),f\>)
 \ar=\ar
G(\<\mu,f\>) + \int_0^t G^\prime(\<Z(s),f\>)\<Z(s),Af\> ds \cr
 \ar \ar
+ ~\frac{\sigma^2}{2}\int_0^t G^{\prime\prime}(\<Z(s),f\>)f(c)^2 ds
+ \mbox{mart.}
 \eeqlb
\eproposition

\proof Let $\{Z_0(t): t\ge 0\}$ be the process mentioned in
Lemma~\ref{l4.7} and let $Z(t) = P_t\mu + Z_0(t)$. Then (\ref{4.10})
clearly holds. \qed

\bproposition\label{p4.9} Let $\{Z(t): t\ge 0\}$ be a solution to
the martingale problem (\ref{4.10}) with sample paths in
$C([0,\infty)$, $\mcr{S}^\prime(\mbb{R}))$. Then we have the
Langevin type stochastic equation
 \beqlb\label{4.11}
\<Z(t),f\> = \<\mu,f\> + \sigma B(t)f(c) + \int_0^t \<Z(s),Af\> ds,
\quad t\ge 0, f\in \mcr{S}(\mbb{R}),
 \eeqlb
where $\{B_t: t\ge 0\}$ is a standard one-dimensional Brownian
motion. \eproposition

\proof By applying (\ref{4.10}) to suitable truncations of the
function $G(z)=z$ we get
 \beqlb\label{4.12}
\<Z(t),f\> = \<Z(0),f\> + M_t(f) + \int_0^t \<Z(s),Af\> ds,
 \eeqlb
where $\{M_t(f)\}$ is a local martingale. By It\^o's formula,
 \beqnn
\<Z(t),f\>^2 = \<Z(0),f\>^2 + 2\int_0^t \<Z(s),f\>\<Z(s),Af\> ds +
\<M(f)\>_t + \mbox{local mart.}
 \eeqnn
On the other hand, if we apply (\ref{4.10}) directly to suitable
truncations of the function $G(z)=z^2$, then
 \beqnn
\<Z(t),f\>^2 = \<Z(0),f\>^2 + 2\int_0^t \<Z(s),f\>\<Z(s),Af\> ds +
\sigma^2f(c)^2t + \mbox{local mart.}
 \eeqnn
Comparing the above two equations we have
 \beqlb\label{4.13}
\<M(f)\>_t = \sigma^2f(c)^2t, \qquad t\ge 0, f\in \mcr{S}(\mbb{R}).
 \eeqlb
Clearly, (\ref{4.12}) and (\ref{4.13}) determine a continuous
orthogonal martingale measure on $[0,\infty)\times \mbb{R}$ with
intensity $\sigma^2\delta_c(x)dsdx$. By El Karoui and M\'el\'eard
(1990, Proposition~II-1) we have
 \beqnn
M_t(f) = \sigma B_tf(c), \qquad t\ge 0, f\in \mcr{S}(\mbb{R}).
 \eeqnn
for a standard one-dimensional Brownian motion $\{B(t): t\ge 0\}$.
\qed

\bproposition\label{p4.10} Let $\{Z(t): t\ge 0\}$ be a solution to
the stochastic equation (\ref{4.11}) with sample paths in
$C([0,\infty)$, $\mcr{S}^\prime(\mbb{R}))$. Then we have a.s.\
 \beqlb\label{4.14}
\<Z(t),f\> = \<\mu,P_tf\> + \sigma\int_0^t P_{t-s}f(c) dB(s), \quad
t\ge 0, f\in \mcr{S}(\mbb{R}).
 \eeqlb
\eproposition

\proof If $\{Z(t): t\ge 0\}$ is a solution of (\ref{4.11}) with
sample paths in $C([0,\infty)$, $\mcr{S}^\prime(\mbb{R}))$, we have
 \beqnn
\int_0^t\<Z(s),P_{t-s}f\> ds
 \ar=\ar
\int_0^t\<\mu,P_{t-s}f\> ds + \sigma\int_0^tP_{t-s}f(c)B(s) ds\cr
 \ar \ar
+\int_0^t ds\int_0^s\<Z(u),P_{t-s}Af\> du \cr
 \ar=\ar
\int_0^t\<\mu,P_{t-s}f\> ds + \sigma\int_0^tP_{t-s}f(c)B(s) ds\cr
 \ar \ar
+\int_0^t du\int_u^t\<Z(u),P_{t-s}Af\> ds \cr
 \ar=\ar
\int_0^t\<\mu,P_{t-s}f\> ds + \sigma\int_0^tP_{t-s}f(c)B(s) ds\cr
 \ar \ar
- \int_0^t \<Z(u),f\> du + \int_0^t \<Z(u),P_{t-u}f\> du.
 \eeqnn
It follows that
 $$
\int_0^t\<Z(s),f\> ds
 =
\int_0^t\<\mu,P_{t-s}f\> ds + \sigma\int_0^tP_{t-s}f(c)B(s) ds.
 $$
Consequently, we have
 \beqnn
\int_0^t\<Z(s),Af\> ds
 \ar=\ar
\int_0^t\<\mu,P_{t-s}Af\> ds + \sigma\int_0^tP_{t-s}Af(c)B(s) ds\cr
 \ar=\ar
\<\mu,P_tf\> - \<\mu,f\> - \sigma f(c)B(t) + \sigma\int_0^t
P_{t-s}f(c) dB(s),
 \eeqnn
where in the last equality we have used the formula of integration
by parts. Then we use (\ref{4.11}) again to see (\ref{4.14}) holds.
\qed

A combination of the above propositions shows that the Langevin type
equation (\ref{4.11}) has a pathwise unique solution and the
martingale problem (\ref{4.10}) is well-posed. Moreover, by
(\ref{4.14}) it is easy to show that $\{Z(t): t\ge 0\}$ is a Markov
process with transition semigroup $(Q^c_t)_{t\ge 0}$ defined by
 \beqlb\label{4.15}
\int_{\mcr{S}^\prime(\mbb{R}^d)}\e^{i\<\nu,f\>}Q^c_t(\mu,\d\nu)
 =
\exp\bigg\{i\<\mu,P_tf\> - \frac{\sigma^2}{2}\int_0^t P_sf(c)^2 d
s\bigg\}.
 \eeqlb
A distribution-valued Markov process with transition semigroup in
this form is usually called an Ornstein-Uhlenbeck type process. The
process $\{Z(t): t\ge 0\}$ describes the asymptotic fluctuations of
the single point catalytic SBM as the branching mechanisms are
small. More precisely, we have the following theorem.

\btheorem\label{t4.11} As $k\to\infty$, the process $\{Z_k(t): t\ge
0\}$ converges weakly in $C([0,\infty)$, $\mcr{S}^\prime(\mbb{R}))$
to the Ornstein-Uhlenbeck process $\{Z(t): t\ge 0\}$ with transition
semigroup $(Q^c_t)_{t\ge 0}$ and $Z(0)=0$. \etheorem

\proof By Lemma~\ref{l4.6} the family $\{Z_k(t): t\ge 0\}_{k\ge 1}$
is tight in the space $C([0,\infty), \mcr{S}^\prime(\mbb{R}))$. Then
the result follows from Lemma~\ref{l4.7} and the well-posedness of
the martingale problem. \qed

\end{document}